\documentclass{amsart}
\newcommand{\subsubsect}[1]{\bigskip \noindent {\bfseries #1. }}

\usepackage[frame,cmtip,arrow,matrix,line,graph]{xy}
\usepackage{graphicx}              

\newtheorem{thm}{Theorem}[section]

\theoremstyle{definition}
\newtheorem{Def}[thm]{Definition}

\newtheorem{ex}[thm]{Example}
\theoremstyle{remark}
\newtheorem{rem}[thm]{Remark}

\numberwithin{equation}{section}

\newcommand{\x}{\times}
\newcommand{\op}{\oplus}
\newcommand{\del}{\partial}
\newcommand{\minus}{\backslash}

\newcommand{\CC}{\mathcal{C}}
\newcommand{\F}{\mathbb{F}}

\newcommand{\FI}{\textrm{FI}}

\newcommand{\GG}{\mathcal{G}}
\newcommand{\M}{\mathcal{M}}
\newcommand{\N}{\mathbb{N}}
\newcommand{\Om}{\Omega}

\newcommand{\s}{\sigma}
\newcommand{\Si}{\Sigma}
\newcommand{\Z}{\mathbb{Z}}

\newcommand{\GL}{\operatorname{GL}}

\newcommand{\Aut}{\operatorname{Aut}}
\newcommand{\Diff}{\operatorname{Diff}}

\newcommand{\Stab}{\operatorname{Stab}}
\newcommand{\id}{\textrm{id}}

\newcommand{\inc}{\hookrightarrow}

\newcommand{\rar}{\longrightarrow}

\newcommand{\note}[1]
{\textcolor{blue}{ [#1]}}

\numberwithin{equation}{section}

\begin{document}


\title{Homological stability: a tool for computations}


\author{Nathalie Wahl}

\address{Dept. of Mathematics, University of Copenhagen,
  Universitetsparken 5, 2100 Copenhagen, Denmark} 


\begin{abstract}
Homological stability has shown itself to be a powerful tool for the computation of  homology of families of groups such as general linear groups, mapping class groups or automorphisms of free groups. We survey here tools and techniques for proving homological stability theorems and for computing the stable homology, and illustrate the method through the computation of the homology of Higman-Thompson groups. 
\end{abstract}

\maketitle


\section{Introduction}

Homology is an invariant that comes in many flavours. We will here mostly be concerned with group homology, but the story we will tell can be told in other contexts as well.
Like many invariants, while easy to define, homology is often difficult to compute. What homological stability has shown to us over the years is that in some situations, it is easier to compute infinitely many homology groups at once than computing a single one by itself. We will in this paper illustrate this through examples, and try to give the reader a sense of how to do homological computations using stability methods, and a sense of when such methods are likely to work.

\smallskip

Many mathematical objects come in families. We will here be interested
in families of groups like the symmetric groups $\Sigma_n$, the braid group $B_n$, the general linear groups $\GL_n(R)$ over a ring $R$, or 
automorphism groups $\Aut(F_n)$ of free groups. 
These examples are more than collections of groups: they all have an additional structure in the form of maps
\begin{align*}
 \op\colon& \Si_n\x \Si_m\rar \Si_{n+m} \\
 \op\colon&    B_n\x B_m\rar B_{n+m} \\
 \op\colon&  \GL_n(R)\x \GL_m(R)\rar \GL_{n+m}(R) \\
 \op\colon& \Aut(F_n)\x \Aut(F_m) \rar \Aut(F_{n+m})  
\end{align*}
by ``block sum'' of permutations, braids or matrices, or juxtaposition of  automorphisms. Another important flavour of example for us will be families of  mapping class groups of surfaces or 3-manifolds with sum $\op$ an appropriate boundary connected sum.

Taking $m=1$ in the above and evaluating the maps $\op$ at the identity element in $\Si_1, B_1, \GL_1(R)$ or $\Aut(F_1)$, gives sequences of groups 
\begin{align*}
  & \Si_1\rar \Si_2\rar \Si_3\rar \dots \\
   & B_1\rar B_2\rar B_3\rar \dots \\
  & \GL_1(R)\rar \GL_2(R)\rar \GL_3(R)\rar \dots\\
  &  \Aut(F_1)\rar \Aut(F_2)\rar \Aut(F_3)\rar \dots
\end{align*}
We are here interested in the following property of such sequences of groups: 

\begin{Def}\label{def:stab}
  A sequence of groups $G_1\to G_2\to \dots$ satisfies {\em homological stability} if the associated sequence of homology groups
  \begin{equation}\label{equ:stab}
  H_i(G_1)\to H_i(G_2)\to H_i(G_3)\to \dots
\end{equation}
is eventually constant for each $i$, that is if $H_i(G_n)\xrightarrow{\cong} H_i(G_{n+1})$ for $n$ large enough with respect to $i$. 
\end{Def}
Unless explicitly otherwise stated, homology here means homology with integral coefficients: $H_*(-)=H_*(-;\Z)$.
Typical stability bounds are linear,  of the form $n\ge ki+a$, for $k$ the {\em slope} of stability.  

Definition~\ref{def:stab} clearly makes sense in other contexts, with
the groups and group homology  replaced by some other type of object
and associated homology theory. Much of what we will present here is known to 
generalize to sequences of spaces, and, to some level, sequences of
algebras. We will focus here on the  case of groups for simplicity,
and because it is already very rich.

\smallskip

All the examples of families of groups mentioned above are known to satisfy homological stability. 
  In the 60s, Nakaoka computed the homology of the symmetric groups $\Si_n$ in \cite{Nak60} and observed that 
    $$H_i(\Si_n)\xrightarrow{\cong} H_i(\Si_{n+1}) \ \ \ \ \textrm{for}\ \ i\le \frac{n}{2}.$$
     Arnold computed in \cite{Arn70} that the same holds for the homology of the braid groups and 
    around the same time,  Quillen, interested in algebraic K-theory  \cite{QuiNotes} (see also \cite{SprWahQ}), showed for example that, for a field $\F\neq \F_2$, 
    $$H_i(\GL_n(\F))\xrightarrow{\cong}  H_i(\GL_{n+1}(\F)) \ \ \ \ \textrm{for}\ \  i\le n.$$
 Harer showed in the 80's that also mapping class groups of surfaces satisfy homological stability \cite{Har85}, a result that was extended to non-orientable surfaces by the author \cite{Wah08}. 
    For the automorphisms of free groups, the first proof goes back to Hatcher \cite{Hat95}, while Hatcher and the author proved a very general stability theorem for mapping class groups of 3-manifolds \cite{HatWah10}: if $M,N$ are any orientable $3$--manifolds such that $\del M\neq \emptyset$, then the map $\pi_0\Diff(M\# N^{\# n})\to \pi_0\Diff(M\# N^{\# n+1})$ extending diffeomorphisms by the identity on the added summand $N$, induces an isomorphism on $H_i$ for $i\le \frac{n-2}{2}$. And many more stability results for families of groups are known!

\subsubsect{Quillen's stability argument} 
    Quillen  devised a strategy for proving homological stability 
    using a spectral sequence associated to the action of the groups on appropriate spaces: 
    To apply Quillen's strategy to a  family of groups $\{G_n\}_{n\ge 0}$,  one needs to find a family of simplicial objects $\{W_n\}_{n\in \N}$, with $G_n$ acting on $W_n$, satisfying (roughly) the following: 
  \begin{enumerate}
  \item the action is transitive on vertices, and has ``manageable'' sets of orbits of $p$--simplices for every $p$;
  \item  the stabilizer of a $p$--simplex is a previous group in the sequence, e.g., $G_{n-p-1}$;
    \item each $W_n$ is highly connected. 
    \end{enumerate}
    From this data,  one can construct a spectral sequence with $E^1$--term
$$E^1_{p,q}= \bigoplus_{\s_p}H_q(\Stab(\s_p);\Z_\sigma)$$
where the sum runs over representatives of the orbits of $p$--simplices $\s_p$ in $W_n$. 
The spectral sequence together with conditions (1)--(3) and a few minor additional assumptions, allows then for an inductive argument. (See Section~\ref{sec:stab} for some more details.)

Variants and extensions of this strategy have been applied in a variety of contexts. In addition to the examples already mentioned, stability has been shown using this strategy   for $\GL_n(R)$ for many rings $R$ \cite{Maazen,vdK80},  for other classical groups like symplectic groups, orthogonal groups, unitary groups, see eg.~\cite{Cat07,MvdK,Sah,Vog79}, and many other groups. 
The strategy was also adapted to prove homological stability for moduli spaces of manifolds and configuration spaces \cite{RW09,Per16,GalRW18}, or for certain families of algebras \cite{HepIH,BoyHep20}.
So many stability theorems have been proved using this method that it is difficult to mention them all.

\subsubsect{Stable homology}
Let $$G_\infty:=\bigcup_{n=1}^\infty G_n = \underset{n\to \infty}{\operatorname{colim}}\ G_n$$
to be the limit of the sequence of groups. Homological stability can be reformulated as saying that the map $G_n\to G_\infty$ induces an isomorphism 
$$H_i(G_n)\xrightarrow{\cong} H_i(G_{\infty})$$ in an increasing range of degrees $i\le b(n)$ for $b(n)$ a bound depending on $n$.   The limit homology $H_*(G_\infty)$ is called the {\em stable homology}. The power of homological stability comes from the fact that often, the stable homology is easier to compute, because it most often belongs to the world of spectra, where methods of homotopy theory come into play. 
We give here known stable computations for the examples of families of groups described above.  

The Barratt-Priddy-Quillen theorem identifies the stable homology $H_*(\Si_\infty)$ of the symmetric groups with that of the basepoint component $\Omega_0^\infty \mathbb{S}$ of the infinite loop space 
of the sphere spectrum $\mathbb{S}$. Galatius  showed that the same holds for the stable homology of automorphisms of free groups. Combining these results with the best known homological stability ranges gives 

\begin{thm}\cite{Nak60,HatVog98b,BarPri,Gal11}\label{thm:stabSi}  For all $i\le \frac{n}{2}$, \ \ \
  $H_i(\Si_n)\cong H_i(\Omega_0^\infty \mathbb{S}) \cong H_i(\Aut(F_{n+3})),$
\end{thm}

A direct consequence is  that the stable rational homology of $\Aut(F_n)$ is trivial.
The result also gives that the inclusion $\Si_n\inc \Aut(F_n)$ induces a homology isomorphism in the range $i\le \frac{n-3}{2}$,
a fact we only know through the above stable homology computation.

For the braid groups $B_n$, the corresponding result is

\begin{thm}[\cite{Arn70,CLM}]\label{thm:stabBr} For all $i\le \frac{n}{2}$, \ \ \ 
$H_i(B_n)\cong H_i(\Omega_0^2 S^2).$
\end{thm}

F.~Cohen completely computed  homology of the right hand side, see \cite[Paper III, App A]{CLM}, yielding a full computation of the stable homology of the braid groups.

For $\GL_n(R)$, the relevant spectrum  is the $K$-theory spectrum, and here the flow of information has gone the other way around compared to the above examples: In the case where $R=\F_{p^r}$ is a finite field,
the homology $H_*(\GL(\F_{p^r});\F_\ell)$ was completely computed by Quillen for any prime $\ell\neq p$, a computation he used to deduce information about the $K$--theory spectrum \cite{QuillenK}. 
When $\ell=p$, only the stable homology is fully known:

\begin{thm}\cite{QuillenK,QuiNotes,SprWahQ,GKRW-GL}\label{thm:stabCl}\footnote{Note that the paper \cite{GKRW-GL} use a different stability method than Quillen's, see Section~\ref{sec:E_n}.} 
   \ $H_i(\GL_n(\F_{p^r});\F_p)=0$\ \ for all $i\le n+r(p-1)-3$ if $p^r\neq 2$, and for all $i< \frac{2n-2}{3}$ if $p^r=2$. 
\end{thm}
A similar result holds for symplectic, orthogonal and unitary groups, see \cite{SprWahQ,FiePri}.

For mapping class groups of surfaces, the stable homology was computed by Madsen and Weiss, in a breakthrough work that lead to much progress in manifold topology: 
Denoting by $\Sigma_{g,b}$ an orientable surface of genus $g$ with $b$ boundary components, and by 
$S_{h,b}$ a non-orientable surface of genus $h$ with $b$ boundary components, and combining the Madsen--Weiss theorem with the best known ranges for homological stability, as well as the unoriented versions of these theorems,  we have 

\begin{thm}\cite{Bol12,MadWei02,RW09,Wah08}\label{thm:stabGa}
$$\begin{array}{ll}
  H_i(\pi_0\Diff(\Sigma_{g,b}))\cong H_i(\Omega_0^\infty MTSO(2)) &  i\le \frac{2g-2}{3} \\
   H_i(\pi_0\Diff(S_{h,b}))\cong H_i(\Omega_0^\infty MTO(2)) &  i\le \frac{h-3}{3} 
  \end{array}$$
\end{thm}

Here $MTSO(2)$ or $MTO(2)$ are the Thom spectra  of the orthogonal bundle to the universal bundle over the Grassmannian of oriented or non-oriented 2-planes in $\mathbb{R}^\infty$ respectively.
A direct consequence is that the stable rational homology of these groups are the polynomial algebras $\mathbb{Q}[\kappa_1,\kappa_2,\dots]$ and $\mathbb{Q}[\zeta_1,\zeta_2,\dots]$ respectively, where $|\kappa_i|=2i$ and $|\zeta_i|=4i$. In the oriented case, this rational computation had been a conjecture of Mumford. 
This result was generalized to higher dimensional manifolds of even dimension by Galatius and Randal-Williams \cite{GalRW14,GalRW18} and to odd-dimensional handlebodies by Botvinnik and Perlmutter \cite{BotPer,Per18}.  This has since been used to compute e.g.; homotopy groups of the diffeomorphisms of discs, or give a totally new approach to pseudoisotopy theory \cite{Kra20,KupRW20}. 

In Section~\ref{sec:Thompson}, we will explain such a theorem  for
the Higman-Thompson groups, see Theorem~\ref{thm:stabTh}, which computes, as a corollary, the full homology of Thompson's group $V$. 
And we will explain in Section~\ref{sec:stable} why one should not be surprised to see double or infinite loop spaces in the above statements.

\subsubsect{Content of the paper}
In this article, we  want to adress the following questions:
\begin{enumerate}
\item When can one expect that homological stability holds?
\item How does one find appropriate $G_n$--space for  Quillen's stability argument?
  \item How does one compute the stable homology? 
  \end{enumerate}
Let us though make clear from the start that we will of course not give full answers to any of the three questions.

A priori one only needs a sequence of groups $G_1\to G_2\to\dots$
to talk about homological stability. Following the article \cite{RWW} and its generalization \cite{Kra19}, Sections~\ref{sec:monoidal} of the present paper shows that 
having the additional data of a ``block sum'', as exhibited above for the groups $\Si_n, B_n, \GL_n(R)$ or $\Aut(F_n)$, is enough input to run Quillen's argument in the following sense:
Section~\ref{sec:Wn}, we construct  a canonical {\em space of destabilizations} $W_n$  when the sum operation is braided, and
Theorem~\ref{thm:stab}  in Section~\ref{sec:stab} states that homological stability holds whenever these spaces are sufficiently connectivity.
In Section~\ref{sec:twist}, we explain how homological stability with {\em abelian} and {\em polynomial} coefficients automatically also holds, under the same assumption, see Theorem~\ref{thm:stab2} and \ref{thm:stab3}.

In Section~\ref{sec:stable} we will see that the braiding forces
the stable homology, through the ``group completion theorem'',  to be
that of a double loop space, or  an infinite
loop space when the braiding is a symmetry.

In Section~\ref{sec:Thompson}, we will then explain how all these ideas were used in \cite{SzyWah} to show that the homology of Thompson's group $V$ is trivial, via a stability theorem and stable computation for the more general Higman-Thompson groups.

The article  ends with a short section  addressing the wider perspective.

\section{A general framework for Quillen's stability argument}\label{sec:monoidal}

In this section, we describe a framework in which Quillen's strategy for proving homological stability can always be implemented.

\smallskip

Recall that a {\em groupoid} $\GG$ is a category whose morphisms are all invertible. 
A {\em monoidal groupoid} is a groupoid $\GG$ equipped with a sum
$$\oplus: \GG\x \GG \rar \GG$$
that is associative and unital. It is {\em braided} if it is in addition equipped with an isomorphism
$\s_{A,B}: A\oplus B\rar B\oplus A$ for every pair of objects, satisfying the braid identity $\s_{A,B}\s_{A,C}\s_{B,C}=\s_{B,C}\s_{A,C}\s_{A,B}:A\op B\op C\to C\op B\op A$  and 
such that $$\xymatrix{A\oplus B\ar[r]^-{\s_{A,B}}\ar[d]_{f\oplus g}& B\oplus A\ar[d]^{g\oplus f}\\
C\oplus D\ar[r]^-{\s_{C,D}}& D\oplus C}$$
commutes whenever it is defined, see e.g.~\cite{MacLane}.  The groupoid  is {\em symmetric monoidal} if the braiding squares to the identity.

There are many examples of braided and symmetric monoidal groupoids. Standard examples of interest to us are 
the groupoid of sets with disjoint union, the groupoid of $R$--modules with direct sum, the groupoid of groups with free or direct product,
the groupoid of vector spaces equipped with a symplectic or Hermitian form  with the direct sum, or the groupoid of manifolds of a given dimension with an appropriate connected sum operation.
Each of these examples are actually the groupoid of isomorphisms in a braided or symmetric monoidal category. For us only the isomorphisms will play a role.

\begin{figure}
  \centering
  \includegraphics[width=3cm]{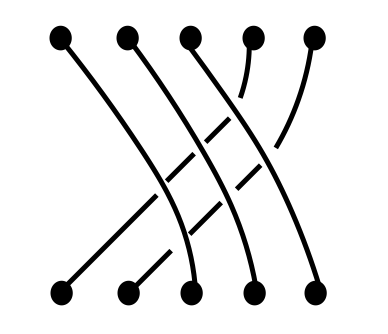}
  \caption{Block braid $b_{3,2}$}\label{fig:braid}
  \end{figure}

\subsubsect{From groups to groupoids}
If we start with a family of groups $\{G_n\}_{n\in \N}$ and defined $\GG=\coprod_nG_n$ to be the groupoid with objects the formal sums $X^{\op n}$ for $n\in \N$ of a generating object $X$, and only non-trivial morphisms $G_n:=\Aut_\GG(X^{\op n})$, then  a monoidal structure on $\GG$, extending the sum in $\N$, is the data of an associative ``sum'' operation   $ G_n\x G_m\xrightarrow{\op}  G_{n+m}$, and a braiding is the data of a homomorphism $\phi:B_n\to G_n$ from the braid group, such that the block braid 
$b_{n,m}$ satisfies that $\phi(b_{n,m})(g\op g')\phi(b_{n,m})^{-1}=g'\op g$ for each $g\in G_n$ and $g'\in G_m$ (see Figure~\ref{fig:braid}).  The groupoid is symmetric precisely if the homomorphism $\phi$ factors through the symmetric group $\Si_n$.

For example, applying this construction to the symmetric groups $\{\Si_n\}_{n\in \N}$ with the block sum of permutations yields the following: the objects are the natural numbers, where we can think of $n$ as representing a set $[n]$ with $n$ elements, and the automorphism group of $[n]$ is $\Si_n$. As $[n+m]\cong [n]\sqcup [m]$, we see that the monoidal sum corresponds to the disjoint union of sets. The resulting groupoid is a skeleton of the groupoid of finite sets of Example~\ref{ex:FB} below. If we instead start with the general linear groups $\{\GL_n(R)\}_{n\in \N}$, we can think of the object $n$ in the resulting groupoid as representing $R^n=R^{\op n}$, whose automorphism group is $\GL_n(R)$. The monoidal product then correspond to the direct sum of $R$--modules, yielding a subcategory of the category of $R$-modules of Example~\ref{ex:RMod} below.

\subsubsect{From groupoids to  groups} If we start instead with a monoidal groupoid $\GG=(\GG,\op)$,  
for any two objects $A,X$ in $\GG$, we get a sequence of groups $G_1\to G_2\to \dots$ by setting 
$$G_n=\Aut_\GG(A\op X^{\op n})$$
with  $$G_n=\Aut_\GG(A\op X^{\op n})\xrightarrow{\_ \op \id_X} G_{n+1}=\Aut_\GG(A\op X^{\op n}\op X).$$
We think of $G_n$ as the automorphism group of ``$A$ stabilized $n$ times in the $X$ direction''.

\begin{ex}\label{ex:FB}
  Let $\GG=\textrm{Sets}^{\textrm{iso}}$ denote the groupoid of finite sets and bijections, with monoidal structure $\op=\sqcup$ given by disjoint union. It is a symmetric monoidal groupoid with the symmetry  the standard bijection $A\sqcup B\xrightarrow{\cong}  B\sqcup A$. Taking $A=\emptyset$ and $X=\{*\}$ yields $G_n=\Si_n$ the symmetric group on $n$ letters, with $\Si_n\to \Si_{n+1}$ the standard inclusion as the subgroup of permutations fixing the last element. 
\end{ex}

\begin{ex}\label{ex:RMod}
  Let $R$ be a ring and let $\GG=R\textrm{--Mod}$ denote the groupoid of $R$--modules and their isomorphisms, with  monoidal product the direct sum $\op$ of modules. This is again a symmetric monoidal groupoid with symmetry given by the standard isomorphism $M\op N\xrightarrow{\cong} N\op M$. 
  Taking $A=0$ and $X=R$, we get $G_n=\GL_n(R)$, the automorphism group of the module $R^{\op n}$, with the map $\GL_n(R)\to \GL_{n+1}(R)$ adding a 1 in the bottom corner of the matrix.
  If we take  $A$ to be any $R$--module, the group $G_n=\GL(A\op R^{\op n})$ is the automorphism group of the module $A$ stabilized $n$ times. 
\end{ex}

\begin{ex}
Let $\GG=\textrm{Groups}^{\textrm{iso}}$ be the groupoid of groups with free product as monoidal structure. This is again a symmetric monoidal groupoid with symmetry the natural isomorphism $G*H\to H*G$. If we take $A=\langle e\rangle$ to be the trivial group and $X=\Z$, we get $G_n=\Aut(F_n)$, the already considered automorphism group of the free group $F_n$. For $A=H$ and $X=G$ any group, the group $G_n=\Aut(H*G*\dots *G)$ is the automorphism group $H$ free product with $n$ copies of $G$, whose stability is studied in \cite{CDG,HatWah10}. 
  \end{ex}

\subsubsect{Modules over monoidal groupoids}
Let  $\GG=(\GG,\op)$ be a monoidal groupoid. 
A category $\mathcal{C}$ is 
a {\em right module over $\GG$} if $\CC$ is equipped with a unital and associative action
$$\mathcal{C}\x \GG \xrightarrow{\op} \CC$$
of  $\GG$. (See \cite[sec 7.1]{Kra19}.) 
Taking $A\in \CC$ and $X\in \GG$ we can again define   $G_n=\Aut_{\CC}(A\op X^{\op n})$, and this  yields just as above a sequence of groups $G_1\to G_2\to \dots$  with the map
$\_\op \id_X: G_n\to G_{n+1}$ adding the identity on the extra copy of $X$.

\smallskip

The sequence of groups $G_n$ obtained above from a monoidal groupoid $\GG$ only is the special case when $\CC=\GG$,  considering $\GG$ as a module over itself. Most of our examples will be of that form, but there are examples from e.g.; manifolds \cite{HatWah10,GalRW18}, or Coxeter groups \cite{Hep16},  that require the more general set-up. (See also \cite{Kra19} for examples in the context of homological stability for topological spaces.)

\subsection{The space of destabilizations}\label{sec:Wn}

Recall from the introduction that to apply Quillen's strategy for proving homological stability, one needs for each $n$ a simplicial object $W_n$ on which the group $G_n$ acts, with appropriate transitivity, stabilizer and connectivity properties. The spaces $W_n$ used in homological stability are most typically one of three types: simplicial complexes, \mbox{(semi-)}simplicial sets, or posets. We will here only discuss spaces of the first two types.  

Ad hoc simplicial objects $W_n$ associated to families of groups $G_n$ have been defined in very many situations to prove stability statements; in fact, most homological stability theorems for families of groups have been so far proved using Quillen's strategy.
Following \cite{RWW}, and its generalization \cite{Kra19}, we construct here the smallest such semi-simplicial set $W_n$  
for any family of groups of the form $G_n=\Aut_\CC(A\op X^{\op n})$ arising as above from the action of a {\em braided} monoidal groupoid $\GG$ on a groupoid $\CC$; the definition of the face maps in $W_n$ will use the braiding of $\GG$. We also define an associated simplicial complex $S_n$.

Fix $\CC$  a module over a braided monoidal groupoid $\GG$, with $A$ an object of $\CC$, and $X$ an object of $\GG$ as above.

\begin{Def}(\cite[Def 2.1]{RWW},\cite[Def 7.5]{Kra19})\label{def:WS}
  The {\em space of destabilizations} $W_n(A,X)_\bullet$ is the semi-simplicial set with set of $p$-simplices 
\begin{align*}
  W_n(A,X)_p 
                 =& \ \{(B,f)\ |\ B\in \textrm{Ob}(\CC) \ \textrm{and}\ f:B\op X^{\op p+1}\to A\op X^{\op n} \ \textrm{in}\ \CC\}/_\sim
\end{align*}
where $(B,f)\sim (B',f')$ if there exists an isomorphism $g:B\to B'$ in $\CC$ satisfying that $f=f'\circ (g\op \id_{X^{\op p+1}})$. 
The face map $d_i: W_n(A,X)_p\to W_n(A,X)_{p-1}$ is defined by $d_i[B,f]=[B\op X, d_if]$ for  
$$d_i f\colon B\op X\op X^{p} \xrightarrow{\id_B\op b_{X^{\op i},X}^{-1}\op \id_{X^{\op p-i}}} B\op X^{\op i}\op X\op X^{\op p-i}\xrightarrow{\ f\ } A\op X^{\op n},$$
for $b_{X^{\op i},X}^{-1}: X\op X^{\op i}\to  X^{\op i}\op X$ coming from the braiding in $\GG$. 
\end{Def}

The group $G_n=\Aut_{\CC}(A\op X^{\op n})$  
acts on $W_n(A,X)_\bullet$ 
by postcomposition.
The following holds for the action: 
\begin{multline}\label{equ:cancel} \textrm{(local cancellation) If }\  
  Y\op X^{\op p+1}\cong  A \op X^{\op n} \ \ \Longrightarrow\ \ Y\cong A\op X^{\op n-p-1}, \\
  \textrm{then }  G_n \textrm{ acts transitively on } W_n(A,X)_p.
\end{multline}
\begin{multline}\label{equ:inject} \textrm{(injectivity) If the stabilization } G_{n-p-1}\to G_{n} \textrm{ taking } f \textrm{ to } f\op \id_{X^{p+1}} \textrm{ is injective,}\\
  \textrm{then there is an isomorphism } G_{n-p-1} \cong \Stab(\sigma_p) \textrm{ for of any } p\textrm{--simplex } \sigma_p.
  \end{multline}
  A direct consequence is that, under these two mild conditions on the $\GG$--module $\CC$, the set of $p$--simplices $W_n(A,X)_p$ of the space of destabilizations is isomorphic to $G_n/G_{n-p-1}$. As we will see in Section~\ref{sec:Thompson} in an example, local cancellation can actually be forced by changing the definition of $\CC$ and $\GG$, declaring in particular
  $A\op X^{\op n}$ and $A\op X^{\op m}$ for $n\neq m$ to be non-isomorphic.
  If the second condition is not satisfied, $W_n(A,X)$ needs to be replaced by a semi-simplicial space in Quillen's argument, see \cite[Sec 7.3]{Kra19}.

\smallskip

\begin{rem}
  In the case of a groupoid $\GG=\CC$ acting on itself, the set of $p$--simplices of $W_n(A,X)$ can be interpreted as the set of morphisms from $X^{\op p+1}$ to $A\op X^{\op n}$ in a category $\langle \GG,\GG\rangle$ constructed from the action, see Appendix~\ref{sec:UG}.
   The face maps are then given by precomposition with standard morphisms $X^{\op p}\to X^{\op p+1}$ in that category. 
\end{rem}

From $W_n(A,X)$, one can also define a simplicial complex $S_n(A,X)$ as follows:
\begin{Def}
Let $S_n(A,X)$ be the simplicial complex with the same vertices as $W_n(A,X)$. A set of vertices $\{x_0.\dots,x_p\}$ spans a $p$--simplex in $S_n(A,X)$ if and only if they are the vertices of  a $p$--simplex of $W_n(A,X)$. 
\end{Def}

We will see in Section~\ref{sec:connect} that it is often equivalent, and more convenient, to work with $S_n(A,X)$ for connectivity questions.

  \begin{ex}
    As in Example~\ref{ex:FB}, consider $(\GG,\op)=(\textrm{Sets}^{\textrm{iso}},\sqcup)$  the symmetric monoidal groupoid of finite sets, seen as a module over itself, with 
       $A=\emptyset$ and $X=\{*\}$, giving $G_n=\Si_n$ the symmetric group. 
       A $p$--simplex $[B,f]$ of $W_n(\emptyset,\{*\})$ is determined by the restriction of the bijection $f:B\sqcup [p+1]\to [n]$ to $[p+1]$. So a $p$--simplex of $W_n(\emptyset,\{*\})$  is an ordered tuple of $p+1$ elements of $[n]=\{1,\dots,n\}$.  The $i$th boundary map forgets the $i+1$st element. This semi-simplicial set is known as the {\em complex of injective words}, and it is known to be $(n-2)$--connected \cite{Far78} (see also \cite[Sec 5.1]{RWW}). The simplicial complex $S_n(\emptyset,\{*\})$ has the same vertices as $W_n(\emptyset,\{*\})$, namely the elements of $[n]$, and $p+1$ such elements form a simplex in $S_n(\emptyset,\{*\})$ precisely when there exist an injective word in these letters, i.e.~if they are distinct. Hence $S_n(\emptyset,\{*\})$ identifies with the $(n-1)$--simplex $\Delta^{n-1}$.
    \end{ex}

    \begin{ex}
           Let $(\GG,\op)=(R\textrm{--Mod},\op)$
            be the symmetric monoidal groupoid of $R$--modules acting on itself, with 
            $A=0$ and $X=R$, giving $G_n=\GL_n(R)$ as in Example~\ref{ex:RMod}. A $p$--simplex $[B,f]$ in $W_n(A,X)$, with $f:B\op R^{p+1}\xrightarrow{\cong} R^n$, is determined by the pair
      $(f(B)<R^n,f|_{R^{p+1}}:R^{p+1}\inc R^n)$.
           The simplicial complex $S_n(A,X)$ thus has vertices  pairs $(H,f)$ with $H<R^n$ and $f:R\inc R^n$ so that $R^n=H\op f(R)$, and vertices $((H_0,f_0),\dots,(H_p,f_p))$ form a $p$--simplex if together the maps $(f_0\op \dots\op f_p):R^p\to R^n$ form an injective map with a complement $H$ such that each $H_i=H\op \bigoplus_{j\neq i}f_j(H_j)$. This complex is very closely related to complexes studied by van der Kallen \cite{vdK80} and Charney \cite{Cha84} and is $\frac{n-a}{2}$--connected for $a$ a constant depending on the stable rank of $R$ (see \cite[Lem 5.10]{RWW}).  The fact that simplices are not just split injective homomorphisms, but rather split homomorphisms {\em with a choice of complement} $H$, makes the simplicial complex more intricate to study, but it forces the stabilizer of a $p$--simplex to be exactly $\GL(H)$,  instead of an affine version of the group, which would be the case if complements had not been chosen. 
    \end{ex}

The simplicial complex $S_n(A,X)$ has appeared in the literature in many examples long before it was defined in the above generality. Here are a few additional examples: for the automorphisms of free groups $\Aut(F_n)$, it is essentially the complex of split factorizations of Hatcher and Vogtmann \cite{HatVog98} (see \cite[Sec 5.2.1]{RWW}), for mapping class groups of surfaces with genus stabilization, this identifies with the tethered arc complex of the same authors \cite{HatVog17} (see \cite[Sec 5.6.3]{RWW}), while the poset of simplices of $W_n(A,X)$ in the case of unitary groups already appeared in \cite{MvdK} (see \cite[Sec 5.4]{RWW}).

\subsection{Homological stability}\label{sec:stab}

Let $\CC$ be a module over a braided monoidal groupoid $\GG$ as above, with $A$ and $X$ objects of $\CC$ and $\GG$ respectively. We have so far associated a sequence of groups   $G_n=\Aut_\CC(A\op X^{\op n})$ to this data, together with  a collection of associated $G_n$--spaces $W_n=W_n(A,X)$ and $S_n=S_n(A,X)$. 
We will now use this as an input for  Quillen's strategy for proving homological stability for the groups $G_n$.
It turns out that $W_n$ is best suited for the spectral sequence argument.

The spectral sequence in Quillen's argument is obtained as follows. 
Let   $E_\bullet G_n$ be a free resolution of $\Z$ as a $\Z G_n$--module, and  let $\widetilde{C}_*(W_n)$ denote the augmented cellular complex of $W_n$. Tensoring these two objects together, we get a first quadrant double complex 
$$C_{\bullet,*}=E_\bullet G_n\otimes_{G_n}\widetilde{C}_*(W_n).$$
 Filtering $C_{\bullet,*}$ in the first direction gives a spectral sequence whose $E^1$--page is trivial in a range under the assumption that $W_n$ is highly connected, from which it follows that  the  spectral sequence coming from filtering in the second direction must converge to zero in a range.
 By transitivity of the action, this latter sequence has $E^1$--term
  $$E^1_{p,q}= H_q(\Stab(\s_p))\cong H_q(G_{n-p-1})$$
  under the local cancelation and injectivity assumption of Section~\ref{sec:Wn}, where there are no twisted coefficients because the stabilizer of a $p$--simplex in $W_n$ always fixes the simplex pointwise.  This spectral sequence allows for an inductive argument. 
This argument has been written in full details many places, see \cite[Thm 3.1]{RWW} for the case where $W_n$ is precisely the complex of destabilization considered here, or 
e.g.~\cite[Thm 5.1]{HatWah10} for a version adaptable to more general simplicial objects $W_n$.

    \begin{thm}\cite[Thm 3.1]{RWW}\label{thm:stab}
      Let $G_n=\Aut_{\CC}(A\op X^{\op n})$ for $\CC,\GG,A$ and $X$ as above satisfying  (\ref{equ:cancel}) and  (\ref{equ:inject}),
      and assume that for all $n\ge 0$, there is a $k\ge 2$ such that the space $W_n(A,X)$ is $\frac{n-2}{k}$--connected.
  Then the stabilization map
  $$H_i(G_n)\rar H_i(G_{n+1})$$ is an isomorphism for $i\le \frac{n-1}{k}$ and a surjection for $i\le \frac{n}{k}$.
  \end{thm}

  \begin{rem}
   The paper \cite{RWW} has two additional assumptions on $\GG$: it should have no zero divisors and the unit has no non-trivial isomorphisms, but, as pointed out by Krannich in \cite[Sec 7.3]{Kra19}, these two assumptions are not actually necessary. Indeed,  these assumptions ensure that $\Aut_\GG(A\op X^{\op n})\cong \Aut_{U\GG}(A\op X^{\op n})$ for $U\GG=\langle \GG,\GG\rangle$ a certain category associated to the groupoid $\GG$ (see Section~\ref{sec:UG}), but in fact stability just holds for both groups whether they are equal or not, with the same proof.  
The paper \cite{RWW} also only formulates the result for the case of $\GG$ acting on itself, but the proof generalizes with no significant change, as noted in \cite{Kra19}. 
\end{rem}

\begin{rem}[Stability slope]
The slope $k$ of stability given by the theorem depends on the slope of connectivity of the spaces $W_n(A,X)$, though with the constrain that the best possible slope is slope $2$. This last restriction is due to the structure of the spectral sequence. To obtain a better slope than slope 2 with the spectral sequence described here, one needs additional information about the groups or differentials appearing in the spectral sequence; such better slopes do not follow from a direct inductive argument. 
  \end{rem}

   It is an open question whether stability holds if and only if the spaces $W_n(A,X)$ are highly connected, see \cite[Conj C]{RWW}.

\subsubsect{Connectivity of buildings}\label{sec:connect}
Stability can only be proved using the above argument under the condition that the spaces $W_n$ (the above defined spaces of destabilizations or some other appropriate buildings) are highly connected. This is a place where work that depends on the groups in question comes in. Under mild conditions, the connectivity of $W_n(A,X)$ is controlled by that of $S_n(A,X)$, and $S_n(A,X)$ will also typically be (weakly) Cohen-Macaulay, a very useful property in connectivity arguments, see  \cite[Sec 2.1]{RWW}.

There are a  few general useful facts and tricks that are good to know when working on connectivity questions for such simplicial complexes or semi-simplicial sets, see e.g., \cite[Sec 2]{HatVog17}, \cite[Sec 2,4,5]{EbeRW} or \cite[Sec 3]{HatWah10}, for expositions of tools and techniques. 
For an example of how such arguments look like,  the survey paper \cite{Wah13mcg}  gives a proof of high connectivity of simplicial complexes of arcs relevant for the stability of the mapping class groups of surfaces, assembling tricks and techniques from the literature.

\subsection{Twisted coefficients}\label{sec:twist}

Homological stability is also often considered in the context of homology with twisted coefficients:  Given a sequence of groups $G_1\to G_2\to \dots$, and a sequence of modules $M_1\to M_2\to \dots$ such that $G_n$ acts on $M_n$ and the map $M_n\to M_{n+1}$ is equivariant with respect to the map $G_n\to G_{n+1}$, one can ask whether the resulting sequence
$$H_i(G_1,M_1)\to H_i(G_2,M_2)\to H_i(G_3,M_3)\to \dots$$
stabilizes. We explain briefly here how  the same assumptions as Theorem~\ref{thm:stab} yield that stability also holds for certain types of ``abelian''  and ``polynomial''coefficients.

\subsubsect{Abelian coefficients}
Suppose that $M$ is a $G_\infty$--module. Then we can consider $M$ as a $G_n$--module via the maps $G_n\to G_\infty=\cup_n G_n$. If we write $M_n$ for this module, this gives an example of a compatible family of coefficients for the groups $G_n$. We say that $M$ is {\em abelian} if the action of $G_\infty$ factors through its abelianization $H_1(G_\infty)$.

  \begin{thm}\cite[Thm 3.4]{RWW}\label{thm:stab2}
  Let $G_n=\Aut_\CC(A\op X^{\op n})$ be as in Theorem~\ref{thm:stab} and assume that for all $n\ge 0$, there is a $k\ge 3$ such that the space $W_n(A,X)$ is $\frac{n-2}{k}$--connected.
  Then for any  $H_1(G_\infty)$--module $M$, the stabilization map
  $$H_i(G_n;M)\rar H_i(G_{n+1};M)$$ is an isomorphism for $i\le \frac{n-k}{k}$ and a surjection for $i\le \frac{n-k+2}{k}$.
  \end{thm}

  The simplest example of such an abelian coefficient system is $M=\Z H_1(G_\infty)$. Because untwisted homological stability gives that  $H_1(G_\infty)\cong H_1(G_n)$ for $n$ large enough, we have that  the twisted homology in that case computes the homology of the commutator subgroup. A direct corollary is thus that, under the same hypothesis as Theorem~\ref{thm:stab} (with $k\ge 3$), homological stability also holds for the commutator subgroups $G_n'$: the stabilization map also induces isomorphisms 
  $$H_i(G_n')\xrightarrow{\cong} H_i(G_{n+1}')$$  for $i\le \frac{n-k}{k}$ and a surjection for $i\le \frac{n-k+2}{k}$. This gives for example  homological stability for alternating groups (=commutator subgroups of symmetric groups), or  special automorphism groups of free groups (=commutator subgroups of $\Aut(F_n)$).

 Note that the best possible slope given by the statement is now slope 3.  This is optimal as stated because we know from \cite[Prop B]{Hau78} that slope 3 is optimal for alternating groups, despite the fact that the spaces $W_n(A,X)$ in this case are slope 2 connected.

\subsubsect{Polynomial coefficients} 
Twisted coefficients classically used in homological stability have been of ``polynomial type'', as introduced by Dwyer in \cite{Dwy80} in the case of general linear groups. It turns out that polynomiality in the sense of Dwyer makes sense in our current general framework of groups of the form  $G_n=\Aut_\CC(A\op X^{\op n})$, as we explain now.

To define a coefficient system for the groups $G_n$, we need the data of a module $M_n$ over $G_n$ for each $n$, and a map $M_n\to M_{n+1}$ compatible with the actions. We will here encode this data in a functor from a 
category built from the $\GG$--module $\CC$, in  similar fashion as the spaces $W_n(A,X)$ were build from $\GG$ and $\CC$\footnote{This is again an example of a bracket construction for categories, as described in Section~\ref{sec:UG}}: Let $\CC_{A,X}$ be the category with objects  $A\op X^{\op n}$ and morphisms from  $A\op X^{\op m}$ to  $A\op X^{\op n}$ empty unless   $m\le n$, in which case a morphism is an equivalence class  of maps $f:A\op X^{\op n}\to A\op X^{\op n}$ in $\CC$, with $f\sim f'$ if there is an isomorphism $g:X^{\op n-m} \to X^{\op n-m}$ in $\GG$ such that $f=f'\circ (\id \op g)$.

A functor $M:\CC_{A,X}\to R\textrm{--Mod}$ defines a coefficient system in the above sense, by setting $M_n:=M(A\op X^{\op n})$. 
Because of the equivalence relation in the definition of the morphisms in $\CC_{A,X}$, such coefficient systems have the particularity that $G_m$ acts trivially on the image of the map $M_n\to M_{n+m}$; they are in fact characterized by this property \cite[Prop 4.2]{RWW}.

Using the  braiding of $\GG$,  we can define a functor
$$\Si_X:\CC_{A,X} \to \CC_{A,X}$$ that adds a copy of $X$ ``to the left'', taking $A\op X^{\op n} $ to   $A\op X^{\op n+1}$, and a morphism $f$ to the composition $(\id_A \op b_{X,X^{\op n}}) \circ (f\op \id_X)\circ (\id_A \op b^{-1}_{X,X^{\op n}})$. This functor comes with a natural transformation $\sigma_X\colon \id \Rightarrow \Si_X$ (see \cite[4.2]{RWW}). 
For $M:\CC_{A,X}\to R-\textrm{Mod}$, we define its suspension
$$\Si M=M\circ \Si_X:\CC_{A,X}\to R-\textrm{Mod}.$$
It comes with a natural transformation $M\rightarrow \Si M$ induced by $\sigma_X$.

A {\em finite degree coefficient system} is defined inductively as follows: the trivial coefficient system $M\equiv 0$ is by definition of degree $-1$, and a coefficient system $M$ is of {\em degree $r$} if the natural transformation $M\to \Si M$ has trivial kernel, and cokernel of degree $r-1$ \cite[Def 4.10]{RWW}. 
For example, constant coefficient systems are of degree $0$,  and all
finitely presented $FI$--modules are coefficient systems of finite
degree for the symmetric groups \cite[Prop 4.18]{RWW}. The Burau representation of the braid group is an example of a coefficient system of degree~1 \cite[Ex 4.15]{RWW}.

  \begin{thm}\cite[Thm A]{RWW}\label{thm:stab3}
 Under the same hypothesis as Theorem~\ref{thm:stab},  if $\{M_n\}_{n\in \N}$ is a polynomial coefficient system of degree $r$, then
   $$H_i(G_n;M_n)\rar H_i(G_{n+1};M_{n+1})$$ is an isomorphism for $i\le \frac{n}{k}-r-1$ and a surjection for $i\le \frac{n}{k}-r$.
  \end{thm}

\section{Group completion and the stable homology}\label{sec:stable}

The fact that the groupoid $\GG$ is braided or symmetric monoidal has direct implications for the stable homology of the groups $G_n=\Aut_\CC(A\op X^{\op n})$ we have been considering here.  
We briefly discuss here the case of automorphism groups $G_n=\Aut_\GG(X^{\op n})$ in $\GG$, and refer to  \cite[Sec 3.2]{RWW} for more details, and for some words about the case $G_n=\Aut_\GG(A\op X^{\op n})$.

\subsubsect{$E_n$--algebras}\label{sec:E_n}
A (topological) $E_n$--algebra is an algebra over the little $n$--disc operad. When $n=1$, such an object goes also under the name $A_\infty$--algebra; it is a space with a multiplication that is associative ``up to all higher homotopies''. When $n\ge 2$, the multiplication is in addition homotopy commutative, with ``more and more'' homotopies as $n$ grows, all the way to an $E_\infty$--algebra that is  commutative up to all higher homotopies. In particular, any $E_n$--algebra is a topological monoid, that is homotopy commutative whenever $n\ge 2$.
(See e.g., \cite{BoaVog,May72}.)

These algebraic structures are relevant for us for the following reason: 
the geometric realization  $|\mathcal G|$ of the nerve  of a monoidal, braided monoidal, or symmetric monoidal category $\mathcal G$
is respectively an $E_1$--, $E_2$-- or $E_\infty$--algebra, see eg., \cite{May74}, \cite[Sec 8]{FSV}. When $\CC$ is a module over a braided monoidal groupoid $\GG$, then $|\CC|$ is an {\em $E_1$--module over the $E_2$--algebra $|\GG|$} in the sense of \cite{Kra19}.

\smallskip

The primary example of an $E_n$--algebra is the $n$--fold loop space $\Om^n X=\textrm{Maps}_*(S^n,X)$ of a space $X$. For $n=\infty$, an $\infty$--loop space is an $n$--fold loop space $Y=\Om^nX_n$ for every $n$, where the spaces $X_n$ together form a {\em spectrum} $\mathbb{X}$.  
Loops have the particularity that they possess homotopy inverses with respect to concatenation, which is the monoid structure underlying their $E_n$--algebra structure. The {\em recognition principle} for iterated loop spaces says that, after ``group completion'', i.e.~after adding homotopy inverses, any $E_n$--algebra is an $n$--fold loop space \cite{May72} (see also \cite{BoaVog,Seg73}).
  Explicitly, the group completion of a topological monoid $(M,\oplus)$ is the space $\Om B_\op M$, where $B_\op$ denotes the bar construction, a simplicial space constructed from $M$ and the sum $\op$. The {\em group completion theorem} states that, if $(M,\op)$ is homotopy commutative, then $H_*(\Om B_\op M)\cong H_*(M_\infty)$ for $M_\infty$ an appropriate ``limit'' space defined from $M$, see \cite{McDSeg,RW13}.

Applying  this to the realization $|\GG|$ of a braided  monoidal groupoid, we get that its group completion $\Om B_\op|\GG|$ is 
a double loop space $\Om^2 X$, or an infinite loop space $\Om^\infty \mathbb{X}$ if $\GG$ was actually symmetric. For $\GG$ of the form $\GG=\coprod_{n\ge 0}G_n$  with $G_n=\Aut_{\GG}(X^{\op n})$, the limit space $|\GG|_\infty$ identifies with $\Z\x BG_\infty$
for $G_\infty=\cup_nG_n=\operatorname{colim} (G_0\to G_1\to G_2\to\dots)$, and 
the group completion theorem thus takes the form
$H_*(\Omega B_{\oplus}|\GG|)\cong  H_*(\Z \x BG_\infty)$. 
Equivalently, it gives that the stable homology of the groups $G_n$ has the following form: 
\begin{align*} H_*(G_\infty)\cong H_*(\Omega_0 B_{\oplus}|\GG|)\cong \left\{\begin{array}{ll}H_*(\Om^2_0X) & \ \textrm{if}\ \GG\  \textrm{is braided} \\
  H_*(\Om^\infty_0\mathbb{X}) &  \ \textrm{if}\ \GG\  \textrm{is symmetric}\end{array}\right.
  \end{align*}
for some space $X$, respectively spectrum $\mathbb{X}$,   
just as in the examples we have seen so far, namely  Theorems~\ref{thm:stabSi}--\ref{thm:stabGa}.
The work in identifying the stable homology of a family of groups  thus comes down, through these classical results, to the question of identifying certain double or, most often, infinite loop spaces arising as classifying spaces of groupoids.
Considered very broadly, this is the subject of $K$--theory. 
In Section~\ref{sec:Thompson}, we sketch one such computation.

\subsubsect{``Higher'' stability and the $E_k$--splitting complex}
The stabilization maps we study here only use a very small part of the $E_2$-- or $E_\infty$--structures we have at hand: taking the sum $\op X$ 
 just uses  part of the underlying $E_1$--module structure. 
The space of destabilizations $W_n(A,X)$ associated to the $E_1$--module $|\CC|$ over the $E_2$--algebra $|\GG|$ and the elements $A\in |\CC|$ and $X\in |\GG|$,   can be thought of as a form of resolution of the space $\coprod_nBG_n$ as $E_1$--module generated by $A$ over the $E_2$--algebra generated by $X$.

One can ask whether there are interesting ``higher'' stabilization maps, summing for example with higher dimensional homology classes, or whether the whole $E_k$-structure can tell us more about the homology of the family of groups. The answer is yes, and is the subject of the body of work \cite{GKRW-Ek,GKRW-GL,GKRW-MCG,GKRW-GLinf} (see also \cite{MilWil} in the context of representation stability).
Considering the full $E_k$--structure has turned out to be powerful, and these papers 
manage to go further than with the classical arguments, including to obtain information about the homology past the stable range. (See also the related paper \cite{Hep20}.) 
The authors define  {\em $E_k$--splitting complexes} that resolve the full $E_k$--structure.
For a relationship between the connectivity of the spaces $W_n(A,X)$ defined here and that of the $E_1$--splitting complex, see  \cite[Thm 13.2]{Hep20}.

\section{Higman-Thomson groups}\label{sec:Thompson}

Sometimes, homological stability is useful in unexpected situations, as turned out to be the case in the study of the homology of Thompson's group $V$.
Thompson's groups come in three flavours: $F<T<V$ where $F$ is a subgroup of the piecewise-linear homeomorphisms of the interval, $T$ a subgroup of those of the circle, and $V$ of the homeomorphisms of the Cantor set. The homology of $F$ and $T$ was computed in the 80's by Brown-Geoghegan and Ghys-Sergiescu in \cite{BroGeo,GhySer}. 
Brown proved a few years later  that the rational homology of Thompson’s group $V$ was trivial, and conjectured that it was also integrally trivial \cite{Bro92}.
Brown's conjecture was proved 25 years later by the author and Szymik in \cite{SzyWah} using  the following  unexpected strategy:
\begin{enumerate}
\item $V=V_{1}$ is part of a family of groups $V_{1}\to V_{2}\to V_{3}\to \dots$ that satisfies homological stability;
\item The homology $H_*(V)$ is entirely stable, i.e.~$H_*(V)\cong H_*(V_{\infty})$
\item The stable homology identifies with that of a trivial infinite loop space. 
\end{enumerate}
In fact, we will see below that each group $V_n$ in the sequence is isomorphic to $V$, but the maps $V_n\to V_{n+1}$ are only isomorphisms after passing to homology.
The strategy works more generally to compute the homology of  the Higman-Thompson groups, so we describe it now in more details in that context.

\smallskip
\begin{figure}
  \centering
  \includegraphics[width=0.3\textwidth]{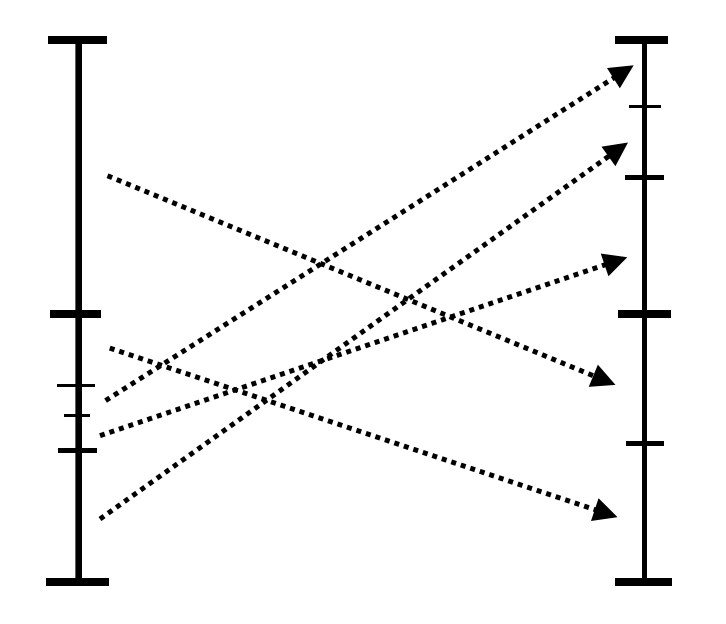}
  \caption{An element of Thompson's group $V=V_{2,1}$ obtained from a binary subdivision of the source and target interval, and a choice of permutation of the subintervals.}\label{fig:V}
\end{figure}

The Higman-Thompson group $V_{k,n}$ is the group of self-maps of a disjoint union of $n$ intervals $I^{\sqcup_n }$ obtained by choosing $k$--ary subdivisions of the source and target, subdividing the interval into $k$ equal sized subintervals and repeating on some of the intervals thus obtained, and matching the resulting subintervals by a chosen bijection. (See \cite[Sec 1.2]{SzyWah}, and Figure~\ref{fig:V} for an example when $k=2$ and $n=1$.)
Thompson's group $V=V_{2,1}$, is the group obtained this way from binary subdivisions of a single interval.  Fixing some $k\ge 2$, we can think of $V_{k,n}$ as the automorphism group of an object $X^{\op n}=I^{\sqcup_n}$ in a groupoid $\mathcal{V}_k=\sqcup_{n\ge 0}V_{k,n}$, just as we have considered in this paper. 
Juxtaposition of intervals induces maps $$V_{k,n}\x V_{k,m}\rar V_{k,n+m}$$
that make the groupoid $\mathcal{V}_k$ symmetric monoidal, with the symmetries coming for block permutations of  the intervals. Hence we can try to apply the stability machine described in the present paper to prove homological stability for the groups $\{V_{k,n}\}_{n\ge 0}$.

Note that there are group isomorphisms $V_{k,n}\cong V_{k,n+(k-1)}$ induced by subdividing an interval into $k$ subintervals, but these isomorphisms are not encoded in the groupoid $\mathcal{V}_k$. For the purpose of homological stability it is convenient to have a rank function, that is, to know what ``$n$'' is at all times. Ignoring these isomorphisms  also gives, by construction, the local cancellation property (\ref{equ:cancel}) which was necessary for the transitivity of the action on the associated complex of destabilization $W_n$.  
So from the point of view of the groups, the objects $I$ and $I^{\sqcup_k }$ are isomorphic, but we will consciously suppress that information in the first part of our argument. 

Let  $W_n=W_n(0,I)$ be the space of destabilizations associated 
to the symmetric monoidal groupoid $\mathcal{V}_k$ (acting on itself) and the objects $0$ and $I$,  and let $S_n=S_n(0,I)$ be its associated simplicial complex, as defined in Section~\ref{sec:Wn}.
The group $V_{k,n}$ can be defined as the automorphism group of an object called the free {\em Cantor algebra} $C_{k}(n)$ of arity $k$ on $n$ generators (see \cite[Def 1.1]{SzyWah}), and 
a $p$--simplex in $W_n$ corresponds to an embedding $C_{k}(p+1)\inc C_{k}(n)$ with complement isomorphic to $C_{k}(n-p-1)$. It is shown in \cite[Cor 3.4]{SzyWah} that $S_n$, and hence also $W_n$ (by \cite[Thm 2.10]{RWW}), is at least $(n-3)$--connected for all $n\ge 2$. The complex $S_n$ has dimension $n-1$ and the idea of the proof of connectivity is to work with its $(n-2)$--skeleton, as simplices that are not maximal correspond to embeddings that have a complement of rank at least 1, i.e.~at least as big as $C_{k}(1)$.  But there are isomorphisms $C_{k}(1)\cong C_{k}(1+(k-1))\cong C_{k}(1+2(k-1))\cong \dots$ so that in practice, a non-trivial complement is actually a  complement that is ``as large as one likes'', which is useful for coning off simplices.

Applying Theorem~\ref{thm:stab} we immediately get that the stabilization map $V_{k,n}\to V_{k,n+1}$ that adds the identity on the new interval, induces an isomorphism $H_i(V_{k,n})\xrightarrow{\cong} H_i(V_{k,n+1})$ in a range increasing with $n$. Coupling this with the fact that the isomorphisms $V_{k,n}\cong V_{k,n+(k-1)}\cong V_{k,n+2(k-1)}\cong \dots$ can be chosen compatibly with the stabilization maps, we get that the rank $n$ can be assumed as large as one like, so that  the  isomorphism $H_i(V_{k,n})\xrightarrow{\cong} H_i(V_{k,n+1})$ actually holds without any bound.

\smallskip

It remains to compute the stable homology. From the results described in Section~\ref{sec:stable}, given that $\mathcal{V}_k$ is a symmetric monoidal groupoid, we know that the stable homology of the groups is that of an infinite loop space. Now here it turns out to be more convenient to do the computation using a different symmetric monoidal groupoid whose group completion also yields the stable homology of the groups $V_{k,n}$, namely the groupoid $\overline{\mathcal{V}}_k$ where we now remember the isomorphism $I\to I^{\sqcup_k}$, or equivalently the isomorphisms of Cantor algebras $C_{k}(n)\cong C_{k}(n+(k-1))$. 
Theorem 5.4 of \cite{SzyWah} says that 
$$H_*(V_{k,\infty})\cong H_*(\Om_0 B_{\op}|\overline{\mathcal{V}}_k|).$$
As  $\overline{\mathcal{V}}_k$ is symmetric monoidal, we again have that its group completion is an infinite loop space and what remains is to find out what the corresponding spectrum is.

So now we are in the world of symmetric monoidal categories, and the idea is simply to find a symmetric monoidal category that is equivalent to  $\overline{\mathcal{V}}_k$ as symmetric monoidal category, and hence group completes to the same infinite loop space, but whose associated spectrum is easier to recognize. Our search was guided by the following observation: the category $\overline{\mathcal{V}}_k$  resembles the category of finite sets and isomorphisms, to which one has declared one extra isomorphism, namely that $[1]$ is now isomorphic to $[k]$,
or, after group completion, $[0]$ is isomorphic to $[k-1]$. 
As already mentioned  in the introduction, the spectrum associated to the category of finite sets (or equivalently to the symmetric groups) is the sphere spectrum $\mathbb{S}$. In homotopy theory, we trivialize by taking cofibers, and the cofiber of the map $\mathbb{S} \xrightarrow{(k-1)} \mathbb{S}$  multiplying by $(k-1)$ is a well-known spectrum $\mathbb{M}_{k-1}$ called the Moore spectrum. Making this idea precise, formulating it on the level of symmetric monoidal categories,  and combining it with the homological stability result described above, led to the following result

\begin{thm}\cite{SzyWah}\label{thm:stabTh}
  There are isomorphisms $$H_*(V_{n,k})\cong H_*(\Om^\infty_0\mathbb{M}_{k-1}).$$
  Specializing to the case $k=2$ yields that $H_*(V)=0$ for $*>0$ as the spectrum $\mathbb{M}_1=\operatorname{cofiber}(\mathbb{S}\xrightarrow{\id}\mathbb{S})$ is trivial. 
  \end{thm}

  Note that the homology of $\Om^\infty_0\mathbb{M}_{k-1}$ for $k\ge 3$ is tractable, and we have many tools available to compute it. For example it is immediate that the rational homology of these groups is trivial, but also that the integral homology is not.  We confirm for instance in \cite[Sec 6]{SzyWah} that  $H_1(V_{k,n})=\Z/2$ for $k$ odd and show that
  the first non-trivial homology group in the $k$ even case is $H_{2p-3}(V_{k,n})=\Z/p$ for $p$ the smallest prime dividing $k-1$.

\section{Perspectives }\label{sec:perspecto}

Many stability results have been proved over the past decades, and one is left to wonder how far homological stability methods can reach. We have highlighted here the idea that braidings seem to be relevant. This is however neither a necessary nor a sufficient condition. We give here some examples that tests the limits of stability, as well as a hint to the wider context homological stability can be considered in.

\medskip \noindent
{\bfseries No braiding = no stability?} Such a statement is not going to ever be literally true, but here are some standard types of examples that are good to have in mind: 
The full braid groups $G_n=B_n$ satisfy homological stability but not the pure braid groups $K_n=\ker(B_n\to \Si_n)$. Likewise the general linear groups $G_n=\GL_n(R)$ satisfy stability for many rings $R$ but not the congruence subgroups $K_n=\GL_n(R,I)=\ker(\GL_n(R)\to \GL_n(I))$. There are in fact many examples of that form with a family of groups $K_n<G_n$ with the groupoid $\GG=\sqcup G_n$ braided monoidal while the groupoid $\mathcal{K}=\sqcup K_n$ is monoidal but not braided, and with the family $G_n$ stabilizing but not the family $K_n$. It turns out that such families $\{K_n\}_{n\ge 0}$ often satisfy instead a form of {\em representation stability} in the sense of \cite{CEF},
see also \cite{Farb,Patzt}.

\subsubsect{Braiding $\not \Rightarrow$ stability}
There are very few examples of braided monoidal categories where we know that homological stability for the associated groups $G_n$ does not hold. One such example, constructed by Patzt  \cite{Pat15}, is the following: consider the category of sets, but using the product $\x$ instead of the disjoint union as monoidal structure. This is a symmetric monoidal groupoid, and if we pick $A=[1]$ and $X=[2]$, we get $G_n=\Si_{2^n}$ is the symmetric group on $2^n$ elements. 
The resulting space $W_n(A,X)$ is however disconnected in this case! And indeed, even though the symmetric groups satisfy homological stability, the stabilization maps in this case do not induce isomorphisms; the induced map on first homology is instead the zero map. 
So existence of a  braiding does not imply stability, which in hindsight is probably not surprising.

There are in addition plenty of examples where we have a braided monoidal groupoid at hand but we don't know that stability holds. For example, the category of $R$-modules over any ring $R$ is symmetric monoidal, but stability for the groups $\GL_n(R)$ is essentially only known under the condition that the ring has finite Bass stable range \cite{vdK80}. But examples of rings for which we know that stability for $\GL_n(R)$ does not hold are surprisingly rare; see \cite{Kup21} for  one example of a ring for which $H_1(\GL_n(R))$ does not stabilize.
For mapping class groups or diffeomorphism groups of manifolds, we essentially know stability in full generality in dimension 2 and 3, but in higher dimension, homological stability for the classifying spaces of diffeomorphism groups is only known for  stabilization by connected sums with certain products $S^p\x S^{q}$ \cite{GalRW18,BotPer}.  Similarly, homological stability for the automorphism groups of vector spaces equipped with a form (symplectic, unitary or orthogonal groups), is mostly known in the particular case of stabilizing with the hyperbolic form, see e.g.~\cite{SprWahQ}. 
In all of these cases, we just do not know the connectivity of the complex of destabilizations.

\subsubsect{Homological stability in other contexts}
We have already mentioned a number of stability results for sequences of spaces. The most classical examples are configuration spaces, going back to the work of McDuff, Segal and F.~Cohen in the 70's \cite{Seg73,McD75,CLM}. In other context, examples seem to be more rare so far, but there is currently a growing interest in stability in the homology of families of algebras, see e.g.~\cite{HepIH,BoyHep20,Robin}, and there exists e.g.~some  results for bounded cohomology of groups \cite{HarMen}. These results are of a very similar flavour as what we have described in the present paper.

\appendix

\section{Adding complements categorically}\label{sec:UG}

The semi-simplicial sets $W_n(A,X)$ of Section~\ref{sec:Wn} and the categories $\CC_{A,X}$ used to  define polynomial coefficients in Section~\ref{sec:twist}, were constructed using equivalence classes of maps in the groupoid $\CC$. Both these constructions are related to a categorical construction,  first considered by Quillen in the context of $K$-theory \cite[p 219]{Gra76}. We recall this construction here and give a few examples. The resulting categories will be natural ``homes'' of the spaces $W_n(A,X)$, and for the polynomial twisted coefficients, which gives some insights.  

\smallskip

Let $\M$ be a category, that is a left module over a monoidal groupoid $(\GG,\oplus)$.  We define  a category $\langle \GG,\M\rangle$ as follows:  
  $\langle\GG,\M\rangle$ has  the same objects as $\M$, and  morphisms from $A$ to $B$ are defined as equivalence classes of pairs $(X,f)$ with $X$ an object of $\GG$ and $f:X\op A\to B$ a morphism of $\M$, where $(X,f)\sim (X',f)$ if there is a commuting diagram
$$\xymatrix{X\op A \ar[d]_{g\op \id}\ar[r]^-f & B \\
X'\op A \ar[ur]_-{f'} & } $$
in $\M$. (If $\CC$ is a right module instead, a category  $\langle \CC,\GG\rangle$ is  defined analogously.)  When $\M$ is a groupoid, as will be the case in our examples, the maps $f$ are isomorphisms and the object $X$ can be thought of as a choice of complement for $A$ inside $B$.

We will  here only consider the case where $\M=\GG$ is a monoidal groupoid acting on itself, and denote by $U\GG=\langle \GG,\GG\rangle$ the resulting category.

\begin{ex}\label{ex:UG}
  Let $(\GG,\op)=(\textrm{Sets}^{\textrm{iso}},\sqcup)$ be the monoidal groupoid of finite sets and bijections of Example~\ref{ex:FB}, with the  monoidal structure induced by disjoint unions. Then $U\GG=\FI$ is the category of finite sets and injections. Indeed, any injection $f:A\inc B$ has, up to isomorphism, a unique complement $X=B\minus f(A)$.
\end{ex}

\begin{ex}
Let $(\GG,\op)=(R\textrm{--Mod},\op)$ be the groupoid of $R$--modules  and isomorphisms of Example~\ref{ex:RMod}, with the  monoidal structure given by direct sum. Then $U\GG$ is closely related to the  category sometimes called $VIC$, with the same objects as $R\textrm{--Mod}$ and with morphisms from $M$ to $N$ given by pairs $(H,f)$ with $f:M\to N$ a split injective homomorphism and $H$ a choice of complement in $N$ of the image: $N=H\op f(M)$ (see e.g.~\cite{PutSam}). 
   \end{ex}

  If the monoidal groupoid $(\GG,\op)$ is braided, one can define a monoidal structure on $U\GG$ as follows: on objects the monoidal structure $\op$  is that of $\GG$, and 
  for $[X,f]$ a morphism from $A$ to $B$ and $[Y,g]$ a morphism from $C$ to $D$, we set
  $$[X,f]\op [Y,g]=[X\op Y,f\op g\circ \id_X\op b_{A,Y}^{-1}\op \id_C]:A\op C \to B\op D$$
where we use the braiding to switch $A$ and $Y$ in $X\op Y\op A\op C$ to be able to apply the morphism $f\op g$. 
The category $U\GG$ is not in general braided (see the next example), though it is symmetric when $(\GG,\op)$ is a symmetric monoidal groupoid, see \cite[Prop 1.8]{RWW}.

\begin{ex}
The braid groups $B_n$ form together a groupoid $\mathcal{B}=\sqcup_nB_n$, that is the free braided monoidal groupoid on one element, where the monoidal structure comes from the juxtaposition of braids. The category $U\mathcal{B}$ can be described in terms of braids with free ends: a morphism from $m$ to $n$ for $m\le n$ in $U\mathcal{B}$ is an equivalence class of braid in $B_n$ where the braid has $n-m$ free ends that can freely pass under, but not over, any other strand, see \cite[Sec 1.2]{RWW}. (It can alternatively be defined in terms of embeddings of punctured discs, see \cite[Sec 5.6.2]{RWW}.) The category $U\mathcal{B}$ is not braided monoidal, but only {\em pre-braided} in the sense of \cite[Def 1.5]{RWW}. 
  \end{ex}

  \begin{rem}
The forgetful map $B_n\to \Si_n$ from the braid groups to the symmetric groups induces a map $ U\mathcal{B}\to FI=U(\textrm{Sets}^{\textrm{iso}})$. Because $\mathcal{B}$ is the free braided monoidal category on one object, it encodes all the structure we have when we picked objects $A$ and $X$ in the groupoids $\CC$ and $\GG$ in Section~\ref{sec:monoidal}.  As pointed out in \cite[Rem 2.8]{Kra19}, the reason we can construct a semi-simplicial set $W_n(A,X)$ comes from the following: 
A semi-simplicial set is a functor $\Delta_{\textrm{inj}}^{\textrm{op}}\to \textrm{Sets}$ for $\Delta_{\textrm{inj}}$ the category of finite ordered sets and ordered injections. One can consider $\Delta_{\textrm{inj}}$ as a subcategory of the category $FI$ of finite sets and injections. Now while the forgetful map $U\mathcal{B}\to FI$ does not admit a splitting, it does admit a partial splitting, in the form of a functor $\Delta_{\textrm{inj}}\to U\mathcal{B}$, and this partial splitting is what rules the semi-simplicial structure of the space of destabilization $W_n(A,X)$. 
\end{rem}


\section*{Acknowledgements}
My interest in homological stability originates in the work of Tillmann and Madsen-Weiss on the moduli space of Riemann surfaces. My thanks goes to them, as well as 
  Ruth Charney, Bill Dwyer, John Harer, Nikolai Ivanov, Dan Quillen and Willem van der Kallen,  for the beautiful papers that inspired much of the research presented here, and of course also  to my stability-collaborators Giovanni Gandini, Allen Hatcher,  Oscar Randal-Williams, David Sprehn, Markus Szymik and Karen Vogtmann.  Finally I would  like to thank S{\o}ren Galatius, Manuel Krannich and Oscar Randal-Williams for thoughtful comments on earlier versions of this paper. 

This work was partially supported by the Danish National Research Foundation through the Copenhagen Centre for Geometry and Topology (DRNF151), and the European Research Council (ERC) under the European Union’s Horizon 2020 research and innovation programme (grant agreement No. 772960). 


\bibliographystyle{plain}
\bibliography{wahlbib}









\end{document}